\documentclass{amsart}
\usepackage{amsmath}
\usepackage{amssymb} 
\usepackage{amscd} 
\usepackage{graphicx}   
\usepackage{xcolor}
\usepackage{marvosym}
\usepackage{float}
\usepackage{ascmac}
\usepackage{wrapfig}
\usepackage{epsfig}
\usepackage{tikz}
\usetikzlibrary{intersections}

\usepackage{hyperref}
\hypersetup{
	colorlinks=true,
	linktoc=all,
    	linkcolor={red!50!black},
   	citecolor={blue!50!black},
  	urlcolor={blue!80!black}
	} 

\usepackage{yhmath} 

\allowdisplaybreaks 

\usepackage{rotating}

\usepackage{booktabs} 

\usepackage{array}
\newcolumntype{M}[1]{>{\centering\arraybackslash}m{#1}} 

\newcommand{\Z}{\mathbb Z}

\newcommand{\bdry}{\partial}
\newcommand{\lk}{\ell k}

\newtheorem{theorem}{Theorem}[section]
\newtheorem{lemma}[theorem]{Lemma}

\newtheorem{proposition}[theorem]{Proposition}
\newtheorem{corollary}[theorem]{Corollary}

\newtheorem{addendum}[theorem]{Addendum}

\newtheorem*{theorem:eta=2}{Theorem~\ref{eta=2}}
\newtheorem*{theorem:wind_stable_wrap}{Theorem~\ref{wind_stable_wrap}}
\newtheorem*{theorem:non_adequate}{Theorem~\ref{non_adequate}}

\theoremstyle{definition}
\newtheorem{definition}[theorem]{Definition}
\newtheorem{example}[theorem]{Example}

\newtheorem{remark}[theorem]{Remark}  

\numberwithin{equation}{section}
\numberwithin{figure}{section}
\numberwithin{table}{section}

\DeclareMathOperator{\wind}{wind}
\DeclareMathOperator{\wrap}{wrap}
\usepackage{fullpage}

\begin{document}

\title{Asymptotic behavior of unknotting numbers \\
of links in a twist family}

\author[K.L. Baker]{Kenneth L. Baker}
\address{Department of Mathematics, University of Miami, 
Coral Gables, FL 33146, USA}
\email{k.baker@math.miami.edu}

\author[Y. Miyazawa]{Yasuyuki Miyazawa}
\address{Department of Mathematical Sciences, Yamaguchi University, 
Yamaguchi 753-8512, Japan}
\email{miyazawa@yamaguchi-u.ac.jp}

\author[K. Motegi]{Kimihiko Motegi}
\address{Department of Mathematics, Nihon University, 
3-25-40 Sakurajosui, Setagaya-ku, 
Tokyo 156--8550, Japan}
\email{motegi.kimihiko@nihon-u.ac.jp}

\begin{abstract}
By twisting a given link $L$ along an unknotted circle $c$, 
we obtain an infinite family of links $\{ L_n \}$.  
We introduce ``stable unknotting number'' which describes the asymptotic behavior of unknotting numbers of links in the twist family.
We show the stable unknotting number for any twist family of links depends only on the winding number  of $L$ about $c$ 
(the minimum geometric intersection number of $L$ with a Seifert surface of $c$) and is independent of the wrapping number  of $L$ about $c$ (the minimum geometric intersection number of $L$ with a disk bounded by $c$).
Thus there are twist families for which the discrepancy between the wrapping number and the stable unknotting number is arbitrarily large. 
\end{abstract}

\maketitle

 \renewcommand{\thefootnote}{}
 \footnotetext{
 }


\section{Introduction}
\label{intro}

For a given link $L$ in $S^3$ with a diagram $D$, 
some set of crossing changes (i.e.\ switching over/under information at a crossing point) transforms $D$ into a diagram of the unlink. 
For each diagram $D$ of $L$, 
we define
\[
u(D) = \mathrm{min}\{ \textrm{the number of crossing changes of}\ D\ \textrm{needed to transform $D$ into a diagram of the unlink}\}.
\]

Using this, the {\em unknotting number} of a link $L$ is defined to be
\[
u(L) = \mathrm{min}\{ u(D) \mid D\ \textrm{is a diagram of}\ L \}.
\]
We may also define the unknotting number of $L$ independently of diagrams 
to be the minimal number of times so that $L$ must be passed through itself to be transformed into an unlink. 
In some literature the value $u(L)$ is called ``unlinking number'', but here we call it ``unknotting number'' even when $L$ is a link to distinguish from ``splitting number'' 
which is the number of crossing changes required to transform $L$ into a link $L'$ for each component $k$ of which we have a $3$--ball containing only $k$. 

\begin{figure}[!ht]
\includegraphics[width=0.45\linewidth]{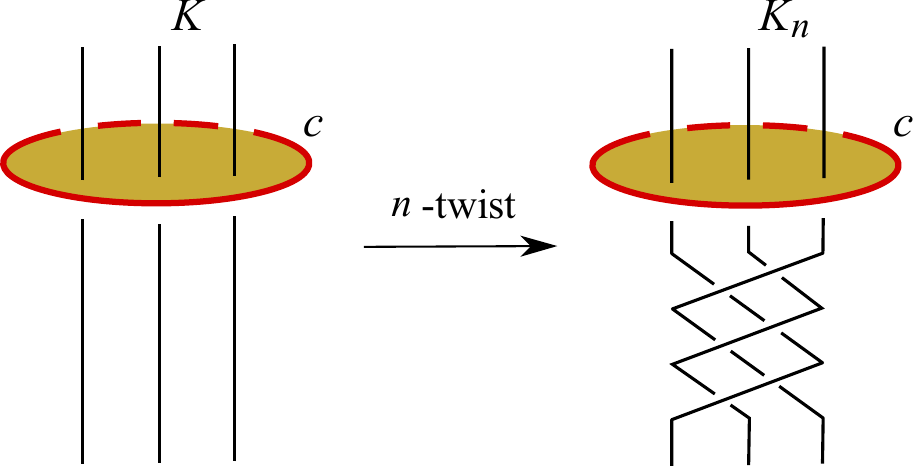}
\caption{Twisting $L$ $n$--times about $c$; $(n = 1)$}
\label{twisting_def}
\end{figure}

Given a link $L$ in $S^3$ and a disjoint unknot $c$, 
then for each integer $n$ a $-1/n$ Dehn surgery on $c$ produces a link $L_{c,n}$ in $S^3$.  
Effectively, $L_{c,n}$ is the result of twisting $L$ about $c$ a full $n$ times.  
This produces a {\em twist family} of links $\{ L_{c, n} \}$ indexed by the integers $n$ where $L=L_{c,0}$; see Figure~\ref{twisting_def}.   
With the twisting circle $c$ understood, we drop it from the notation and speak of the twist family $\{L_n\}$.

For a twist family of links $\{L_n\}$ with twisting circle $c$, the \textit{wrapping number} $\eta = \wrap_c(L)$  and the \textit{winding number} $\omega=\wind_c(L)$ of $L$ about $c$ measure the minimal geometric intersection numbers of $K$ with a disk bounded by $c$ or with any Seifert surface bounded by $c$, respectively.  Hence $0 \leq \omega \leq \eta$ and $\omega \equiv \eta \pmod{2}$.   If $\omega=\eta$, then we say that $L$ is {\em coherent} with respect to $c$ and that the twist family is {\em coherent}.  Note that if $\eta = 0$ or $1$, then $L_{n}=L$ for all $n \in\Z$. Thus we henceforth implicitly assume $\eta \geq 2$.  

\begin{remark}
Usually the winding number is defined as the linking number $\lk(L,c)$ (or its absolute value) when both $L$ and $c$ are oriented.  The above definition corresponds to choosing orientations on $L$ and $c$ so that for each component $L_i$ of $L$ the winding number is non-negative; so that  $\wind_c(L_i) = \lk(L_i,c)  \geq 0$.   
\end{remark}

In this paper we are interested in a behavior of unknotting numbers under twisting operation. 
To investigate an asymptotic behavior of unknotting numbers of links in a twist family, we introduce the stable unknotting number defined below. 

\begin{definition}
Given a twist family of links $\{ L_n \}$, 
let $\{ u(L_n) \}$ be the corresponding sequence of unknotting numbers.
Define the {\em stable unknotting number} of $\{ L_n \}$ to be 
\[u_s(L_n) = \lim_{n \to \infty} \frac{u(L_n)}{n}\]
assuming the limit exists.   
\end{definition}

Now let us take a look at some examples. 

\begin{example}[Torus knots]
Let $V$ be a standardly embedded solid torus in $S^3$, 
and take  a torus knot $T_{p, q}$ on the boundary of $V$ which wraps $q$ times in $V$. 
Let $c$ be an unknot which is a core of the solid torus $S^3 - \mathrm{int}V$. 
Twisting $T_{p, q}$ along $c$, we obtain a twist family of torus knots $\{ T_{p+ qn, q} \}$. 
The Milnor Conjecture \cite{KM} established by Kronheimer and Mrowka asserts that 
the unknotting number of a torus knot $T_{p + qn, q}$ ($p > q \ge 2$) is explicitly determined as $u(T_{p+qn , q}) = \dfrac{|(p+ qn -1)|(q-1)}{2}$.
This enables us to determine the stable unknotting number for a twist family of torus knots. 
\[
u_s( T_{p+ qn, q } ) = \lim_{n \to \infty} \dfrac{|(p+ qn -1)|(q-1)}{2n} = \frac{q(q-1)}{2}.
\]
\end{example}

\begin{example}[Whitehead link and Mazur link]
Figure~\ref{twist_unknotting1} shows two twist families of knots, one generated from the Whitehead link and another from the Mazur link.  
For each of these twist families, 
$u(K_n) \le 1$ for all integers $n$. 
Hence, by definition $u_s( K_n ) = 0$. 

\begin{figure}[!ht]
\includegraphics[width=0.4\linewidth]{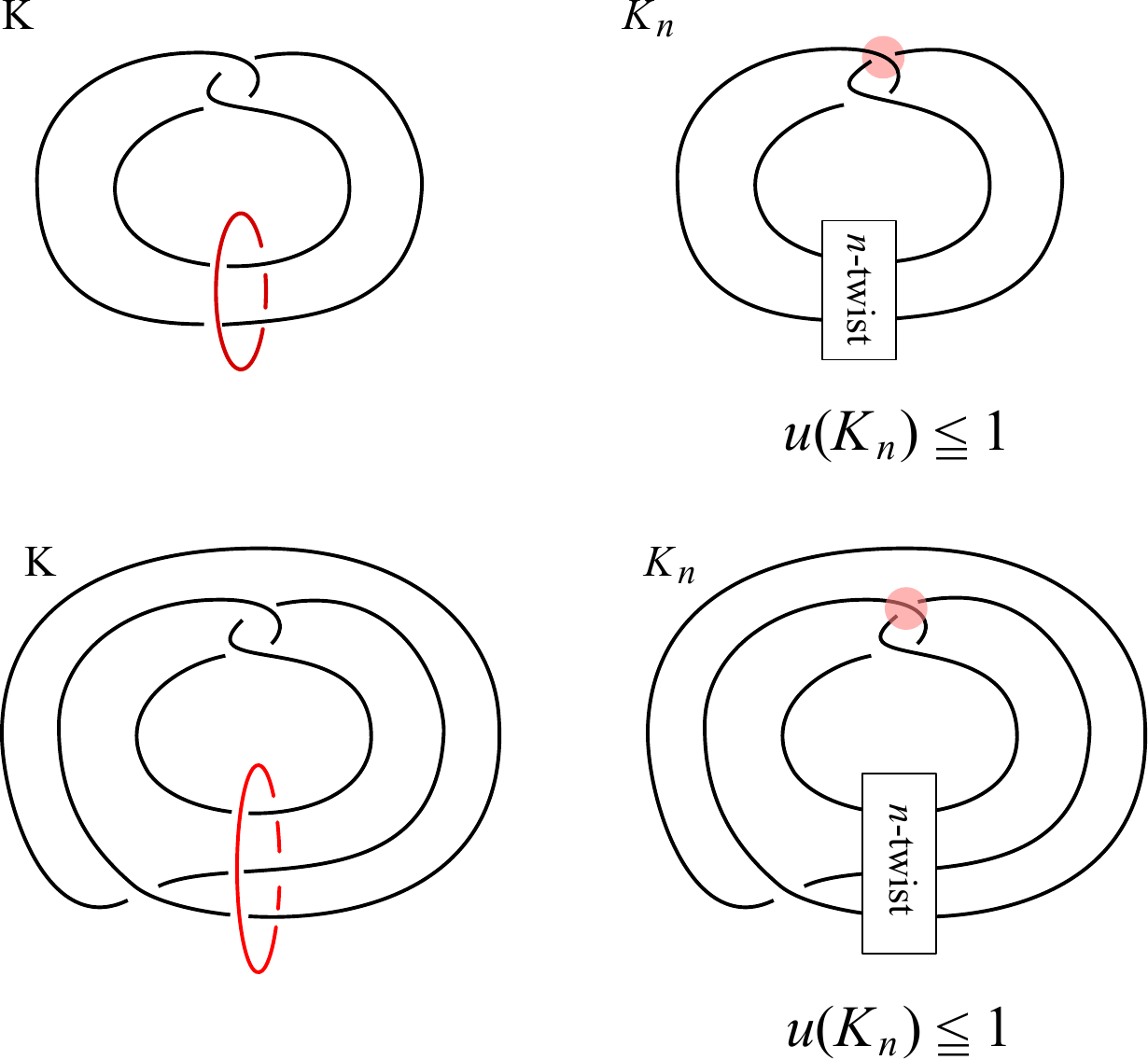}
\caption{Unknotting numbers of the knots $K_n$ is bounded when $K \cup c$ is a Whitehead link (above) or a
Mazur link (below).}
\label{twist_unknotting1}
\end{figure}
\end{example}

Our main result in this paper is the following.

\begin{theorem}
\label{stable_unknotting}
Let $\{ L_n \}$ be a twist family of links with winding number $\omega$. 
Then the stable unknotting number of $\{ L_n \}$ is given by the following. 
\[
u_s( L_n ) = 
\frac{\omega(\omega-1)}{2}.  
\]
\end{theorem}

\begin{addendum}
Let  $\overline{L \cup c} = \overline{L} \cup \overline{c}$ be the mirror image of $L \cup c$, 
and let 
$\overline{L}_n$ be the link obtained from $\overline{L}$ by $n$--twist along $\overline{c}$. 
The winding numbers of $\overline{L}$ about $c$ is also $\omega$. 
Since the mirror image $\overline{L_{-n}}$ of $L_{-n}$ is $\overline{L}_n$, 
whose unknotting number coincides with that of $L_{-n}$. 
Hence, we have 
\[
\lim_{n \to \infty} \frac{u(L_{-n})}{n} = \frac{\omega(\omega-1)}{2}
\]
as well. 
\end{addendum}

\begin{remark}
\label{unknotting_numbers}
The proof of Theorem~\ref{stable_unknotting} shows that 
\[
\frac{\omega(\omega-1)}{2}n + C \le g_4(L_n) \le u(L_n) \le \frac{\omega(\omega-1)}{2}n + C'.
\]
for some constants $C$ and $C'$ that depend  only on the link $L \cup c$. 
Thus  $\displaystyle \lim_{n \to \infty} g_4(L_n)/n = \omega(\omega-1)/2$
which we may regard as a {\em stable slice genus}.  
Hence one may view Theorem~\ref{stable_unknotting} as saying that, for twist families, the stable unknotting number equals the stable slice genus. 
\end{remark}

\begin{remark}
In the proof of Theorem~\ref{stable_unknotting}, Lemma~\ref{u(L_n)-lower} gives the lower bound on $g_4(L_n)$ shown in Remark~\ref{unknotting_numbers}.  This lemma is a direct extension of Proposition~2.4 of \cite{BM} for twist families of knots.
\end{remark}

\begin{corollary}
\label{wind01}
Let $\{ L_n \}$ be a twist family of links with $\wind_c(L) = \omega$. 
The following three conditions are equivalent.
\begin{center}
\begin{tabular}{llllll}
$($1$)$ $u_s( L_n ) = 0$.  & \phantom{SPACE}
$($2$)$ $\omega \le 1$.  & \phantom{SPACE}
$($3$)$ $u(L_n)$ is bounded. 
\end{tabular}
\end{center}
\end{corollary}

\medskip

Using \cite[Theorem]{M} we show that the discrepancy between the wrapping number and the stable unknotting number can be arbitrarily large. 

\begin{proposition}
\label{stable_unknotting_gap}
For given integers $M, N > 0$, 
there exists a twist family of links $\{ K_n \}$ 
with wrapping number $\eta$ and winding number $\omega$
such that 
$\omega^2 \ge M$ and $\eta-\omega^2 \ge N$.   
Hence $\eta-u_s( K_n ) \ge N$ too.
\end{proposition}

\bigskip

\section{Stable unknotting number of knots in a twist family}

In this section we will prove Theorem~\ref{stable_unknotting}. 

Let $L$ be a link in $S^3$. 
Given an unknot $c$ disjoint from $L$, form the link $L \cup c$.  Choose an orientation on $c$.  Then for each component $L_i$ of $L$, choose an orientation so that $\wind_c(L_i) = \omega_i \geq 0$.   With such orientations set, we henceforth regard $L \cup c$ as an oriented link.  Observe that $\wind_c(L) = \omega = \sum \omega_i$  and that twisting $L$ about the unknot $c$ produces the twist family of oriented links $\{L_n\}$ with winding number $\omega$.

\begin{proof}[Proof of Theorem~\ref{stable_unknotting}.]

First we obtain a lower bound on $u(L_n)$.

 Proposition~2.4 of \cite{BM} can be extended from twist families of knots to twist families of links.    The smooth slice genus of the oriented link $L$, $g_4(L)$ is the minimal genus of a connected oriented surface that is smoothly, properly embedded in $B^4$ and bounded by $L$ in $S^3 = \bdry B^4$. The smooth slice genus is a lower bound on unknotting number, that is 
 $g_4(L) \leq u(L) - |L| + 1 \le u(L)$, 
 where $|L|$ denotes the number of components of $L$.   
 
 \begin{lemma}
 \label{u(L_n)-lower}
For a twist family of oriented links $\{L_n\}$ with winding number $\omega$, 
there is a constant $C$ independent of $n$ such that
 \[C+n \omega(\omega-1)/2 \leq g_4(L_n).\]
Hence 
\[C+n \omega(\omega-1)/2 \leq u(L_n).\]
\end{lemma}

\begin{proof}
Following \cite[Proposition 2.4]{BM} we will show that for integers $n \geq m$ we have
\[(n-m)\omega(\omega-1) - 2\eta -2(|L|-1)\leq s(L_n) - s(L_m) \]
where $\omega = \wind_c(L)$, $\eta = \wrap_c(L)$, and $s$ is the Rasmussen $s$-invariant as extended to oriented links by  Beliakova and Wehrli \cite{BW}.  Set $d$ so that $2d = \eta-\omega$.

As in \cite[Proposition 2.2]{BM} we can construct a planar surface $P$ giving a cobordism from $L_n$ to $L_m \sqcup T_{\eta, (n-m)\eta}$ where $T_{\eta, (n-m)\eta}$ has $d$ components running in one direction and $\omega+d$ in the other.  Since $|\bdry P| = \eta + 2|L|$, $-\chi(P) = \eta + 2(|L| - 1)$.
Therefore \cite[Equation (7.1)]{BW} gives
\[ |s(L_m \sqcup T_{\eta, (n-m)\eta}) - s(L_n)| \leq \eta + 2(|L| - 1)\]
and hence (using \cite[Equation (7.2)]{BW}) we have 
\[s(T_{\eta, (n-m)\eta}) -1 - \eta - 2(|L|-1) \leq s(L_n) - s(L_m). 
\tag{$\star$}\]

As in the proof of \cite[Proposition 2.4]{BM}, we continue to have 
\[(n-m)\omega(\omega-1) +1 - \eta \leq s(T_{\eta, (n-m)\eta}). \]
Therefore with $(\star)$ we have
\[(n-m)\omega(\omega-1)  -2\eta - 2(|L|-1) \leq s(L_n) - s(L_m).\]

Putting $m=0$, 
we obtain 
\[
s(L_0) -2\eta - 2(|L| -1) + n\omega(\omega-1) \le s(L_n).
\] 
Let $S$ be a slice surface of $L_n$. 
Remove an open disk from $S$ to obtain a cobordism $\hat{S}$ from $L_n$ to the trivial knot $O$.  
Apply 
\cite[Equation (7.1)]{BW} 
to the cobordism $\hat{S}$ to see that 
$|s(L_n)| = |s(L_n) - s(O)| \le - \chi(\hat{S}) = 2g_4(L_n) + |L_n| - 1= 2g_4(L_n) + |L| - 1$.
Hence we have
\[s(L_0) - 2\eta - 3(|L|-1) + n\omega(\omega-1) \leq 2g_4(L_n).\]
Setting $2C=s(L_0) - 2\eta -3(|L|-1)$ then gives the desired bound.

The last statement follows since slice genus is a lower bound for unknotting number.
\end{proof}

Now we obtain an upper bound on $u(L_n)$.

\begin{lemma}
\label{reduction_coherent}
Let $L \cup c$ be a link with unknotted component $c$. 
Assume that $\wrap_c(L) = \eta$ and $\wind_c(L) = \omega$. 
We may apply $k$ crossing changes on $L$ to obtain $L' \cup c$ so that 
\begin{enumerate}
\item
$\wrap_c(L') = \wind_c(L') = \omega$, and 
\item
the corresponding $k$ crossing changes on $L_n$ convert it into $L'_n$,
i.e.\ the diagram below commutes. 
\end{enumerate}
\end{lemma}

\begin{eqnarray*}
\begin{CD}
	L \cup c @>n-\textrm{twist along}\ c>> L_n \cup c \\
@V{k\ \textrm{crossing changes $\&$ isotopy}}VV
		@VV{k\ \textrm{crossing changes  $\&$ isotopy}}V \\
	L' \cup c @>>n-\textrm{twist along}\ c> L'_n \cup c
\end{CD}
\end{eqnarray*}

\begin{proof}
Winding number determines the homotopy class of an oriented knot in a solid torus.  
Hence any oriented knot in a solid torus is homotopic to a coherent knot, 
one for which its winding number equals its wrapping number.   
Thus $L$ is homotopic in the exterior of $c$ to a link in which each component is coherent, and this link is further homotopic to a coherent link $L'$ where $\wrap_c(L') = \wind_c(L') = \omega$. 
As such homotopies can be realized by a sequence of crossing changes and isotopy, $L$ may be transformed to the link $L'$ 
by some set of $k$ crossing changes in the exterior of $c$.  Figure~\ref{wind_omega} gives one example of the transition from $L \cup c$ to $L' \cup c$.

\begin{figure}[!ht]
\includegraphics[width=1.0\linewidth]{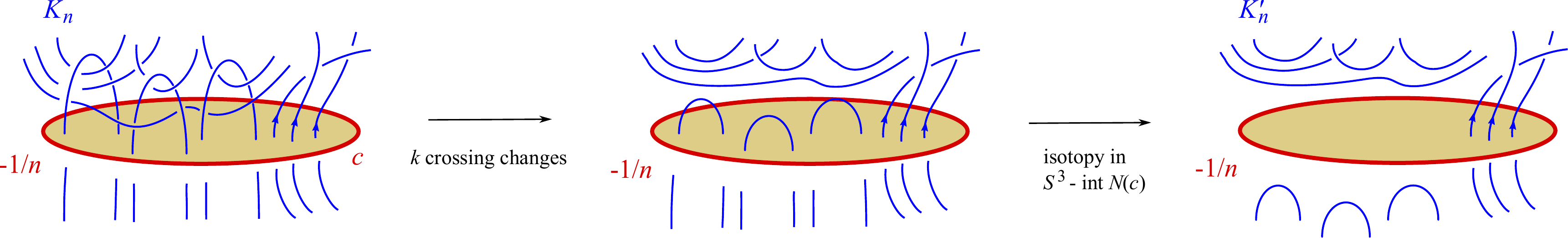}
\caption{$n$--twist and $k$ crossing changes; $\eta = 9, \omega = 3$}
\label{wind_omega}
\end{figure}

Regard $n$--twisting along $c$ as $-1/n$--Dehn surgery along $c$ so that the knots $L_n$ and $L'_n$ are the images of $L$ and $L'$ respectively in $S^3$ after this Dehn surgery.  Thus it follows that the corresponding $k$ crossing changes that transform $L$ to $L'$ in the exterior of $c$ also transform their images $L_n$ to $L'_n$ after $-1/n$--Dehn surgery on $c$.  
\end{proof}

By the construction of the twist family $\{ L'_n \}$ in Lemma~\ref{reduction_coherent}, 
we have both a constant $k \ge 0$ such that 
\begin{equation}
\label{u(L_n)_u(L'_n)}
u(L_n) \le k + u(L'_n)
\end{equation}
for all $n$ and that $\wrap_c(L') = \wind_c(L') = \omega$.

So consider a disk bounded by $c$ that $L'$ intersects $\omega$ times.  Then $L'$ intersects a closed collar neighborhood of this disk (which is homeomorphic to $D^2 \times I$) as a trivial $\omega$-sting braid.  The result of performing $-1/n$ surgery on $c$ then replaces this braid with the braid of $n$ full twists on all $\omega$ of the strands.  Since $\omega(\omega-1)/2$ crossing changes will change the $\omega$-string braid of one full twist back to the trivial braid as demonstrated in Figure~\ref{twist_region_crossing_3}, it follows that
\begin{equation}
\label{u(L'_n)_u(L')}
u(L'_n) \le \frac{\omega(\omega-1)}{2}n + u(L'). 
\end{equation}

\begin{figure}[ht]
\includegraphics[width=0.7\linewidth]{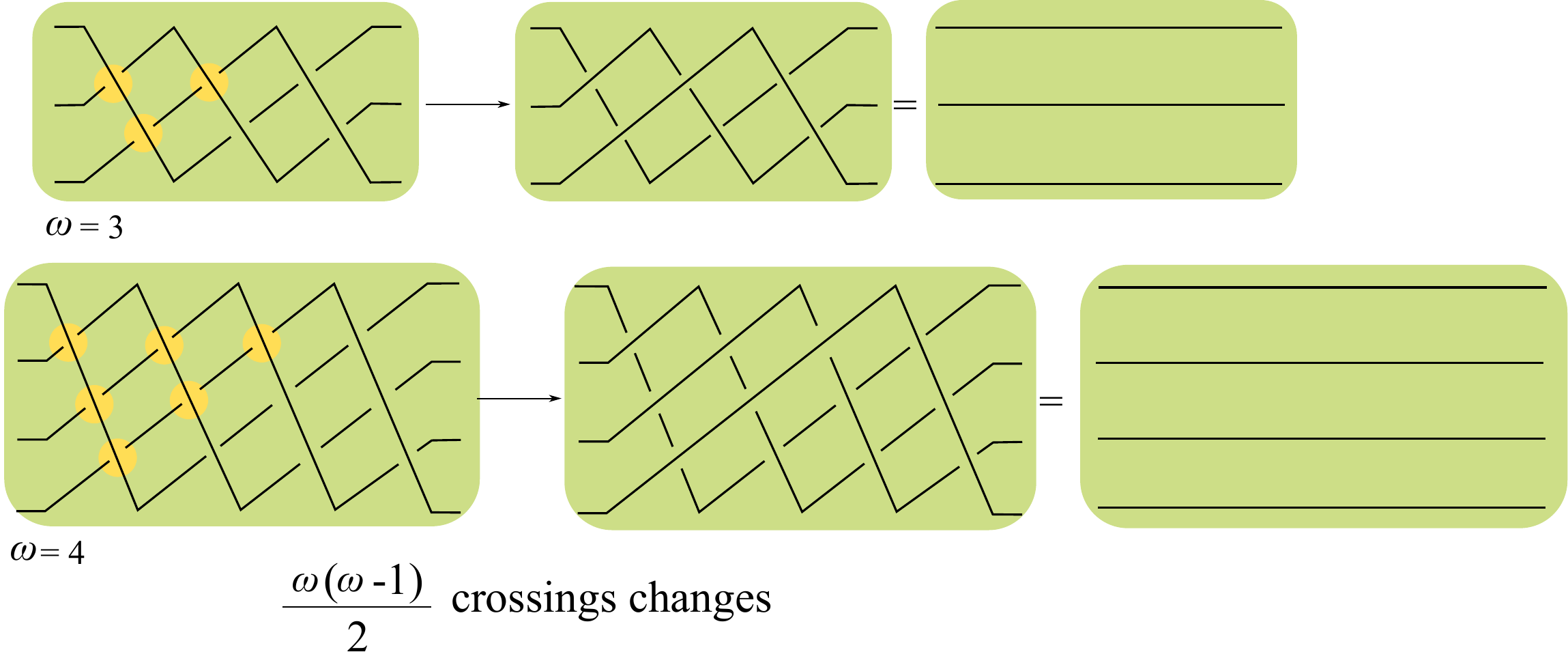}
\caption{ }
\label{twist_region_crossing_3}
\end{figure}

Then combining the inequalities \ref{u(L_n)_u(L'_n)} and \ref{u(L'_n)_u(L')} gives the desired upper bound
\begin{equation}
\label{u(L_n)-upper}
u(L_n) \le k + u(L'_n) \le \frac{\omega(\omega-1)}{2}n + u(L') +k.
\end{equation}

Combining the lower bound in Lemma~\ref{u(L_n)-lower} and the upper bound \ref{u(L_n)-upper} give
\begin{equation}
\frac{C}{2}+  \frac{\omega(\omega-1)}{2}n \le u(L_n)  \le \frac{\omega(\omega-1)}{2}n + u(L') +k
\end{equation}
from which we obtain
\[
\frac{\omega(\omega-1)}{2} 
\le \lim_{n \to \infty}\frac{u(L_{n})}{n} 
\le \frac{\omega(\omega-1)}{2}, 
\]
so that 
\[
u_s( L_n ) = \frac{\omega(\omega-1)}{2}
\] 
as desired.  
\end{proof}

\section{Twist families of links with large wrapping number\\
 and small stable unknotting number}

Although it is quite easy to compute the winding number,  
in spite of simplicity of its definition, it is difficult to determine the wrapping number in general.
When $L$ is not known to be coherent with respect to $c$, one must eliminate the possibility of an isotopy of $L$ reducing the number of intersections of $L$ and a disk bounded by $c$ that apparently realizes the wrapping number. 

Let $\mathcal{L}(\ell, m)$ be the set of links $K \cup c$, where $c$ is a trivial knot, such that 
$\mathrm{wind}_c(K) = \ell$ and $\mathrm{wrap}_c(K) = m$. 
Then following \cite[Theorem]{M} for a given positive integers $\ell, m$ with $\ell \equiv m \pmod 2$
which satisfies (i) $0 = \ell  < m$, (ii) $1 = \ell  < m$, or (iii) $1 < \ell \le m$, 
we may find a link $K \cup c$ in $\mathcal{L}(\ell, m)$. 

\begin{proof}[Proof of Proposition~\ref{stable_unknotting_gap}]
For a given integers $M, N > 0$, 
take a link $K \cup c$ in $\mathcal{L}(\ell, m)$
so that $\ell^2 \ge M$ and $m - \ell^2 \ge N$. 
Since $\mathrm{wind}_c(K) = \ell$ and $\mathrm{wrap}_c(K) = m$, 
the twist family $\{ K_n \}$ satisfies the desired property in Proposition~\ref{stable_unknotting_gap}. 
\end{proof}

\bigskip

\textbf{Acknowledgments}\quad 
KLB was partially supported by the Simons Foundation grant \#523883 and gift \#962034.  He also thanks the University of Pisa for their hospitality where part of this work was done. 

YM has been partially supported by JSPS KAKENHI Grant Number JP22K03315. 

KM has been partially supported by JSPS KAKENHI Grant Number JP19K03502, 21H04428 and Joint Research Grant of Institute of Natural Sciences at Nihon University for 2023. 

We would like to thank the referee for careful reading and useful suggestions.

\end{document}